\documentclass[11pt]{amsart}

\usepackage{enumerate,amssymb}
\usepackage{amsopn,amsfonts,amsmath}
\usepackage[dvips]{graphicx}
\usepackage{float}
\usepackage[american]{babel}
\usepackage{enumerate}
\usepackage{hyperref}
\usepackage{fancyhdr}
\usepackage{color, colortbl}
\usepackage{float}
\allowdisplaybreaks

\def\E{{\mathbb E}}
\def\indic{{\rm {\large 1}\hspace{-2.3pt}{\large
l}}}

\def\R{{\mathbb R}}

\def\xP{{\mathbb P}}

\newtheorem{thrm}{Theorem}[section]
\newtheorem{lmm}[thrm]{Lemma}

\newtheorem{prpstn}[thrm]{Proposition}
\newtheorem{dfntn}[thrm]{Definition}

\begin{document}

\title[]{Pivotal Estimation in High-dimensional Regression via Linear Programming}

\author{Eric Gautier and Alexandre B. Tsybakov}


\maketitle

\begin{center}
CREST (ENSAE),
3 avenue Pierre Larousse, 92 245 Malakoff Cedex, France;\\
\href{mailto:eric.gautier@ensae.fr}{eric.gautier@ensae.fr};
\href{mailto:alexandre.tsybakov@ensae.fr}{alexandre.tsybakov@ensae.fr}.
\end{center}

\begin{abstract}We propose a new method of estimation in high-dimensional linear regression model.
It allows for very weak distributional assumptions including heteroscedasticity, and
does not require the knowledge of the variance of random errors.
The method is based on linear programming only, so that its numerical implementation is faster
than for previously known techniques using conic programs, and it allows one to deal with higher dimensional models.
We provide upper bounds for estimation and prediction errors of the proposed estimator showing that it achieves the same rate as in the more restrictive situation of fixed design and i.i.d.
Gaussian errors with known variance. Following Gautier and Tsybakov (2011), we obtain the results under weaker sensitivity assumptions than the restricted eigenvalue or assimilated conditions.
\end{abstract}

\section{Introduction}\label{s1}
In this paper, we consider the linear regression model
\begin{equation}\label{estruct}
y_i=x_i^T\beta^*+u_i,\ 
\ \quad i=1,\dots,n,
\end{equation}
where $x_i$ are random vectors of explanatory variables in
$\R^p$, and $u_i\in \R$ is a random error. 
The aim is to estimate the vector $\beta^*\in \R^p$
from $n$ independent, not necessarily identically distributed
realizations $(y_i,x_i^T)$, $i=1,\dots,n$.  
We are mainly interested in high-dimensional models
where $p$ can be much larger than $n$ under the sparsity scenario where
only few components $\beta^*_k$ of $\beta^*$ are non-zero ($\beta^*$ is
sparse).

The most studied techniques for high-dimensional regression
under the sparsity scenario are the Lasso, the Dantzig
selector, see, {e.g.}, Cand\`es and Tao~(2007), Bickel, Ritov
and Tsybakov~(2009) (more references can be found in B\"uhlmann
and van de Geer~(2011) and Koltchinskii~(2011)), and
agregation by exponential weighting (see Dalalyan and Tsybakov~(2008), Rigollet and
Tsybakov~(2011, 2012) and the references cited therein).
Most of the literature on high-dimensional regression assumes that the random errors are Gaussian
or subgaussian with known variance (or noise level).  However, quite recently several methods have been proposed
which are independent of the noise level
(see, e.g., St\"adler, B\"uhlmann and van de Geer (2010),
Antoniadis (2010),
Belloni, Chernozhukov and Wang (2011a, 2011b), Gautier and Tsybakov (2011), Sun and Zhang (2011), Belloni, Chen,
Chernozhukov, and  Hansen (2012) and Dalalyan (2012)). Among these, the methods of
Belloni, Chernozhukov and Wang (2011b), Belloni, Chen,
Chernozhukov, and  Hansen (2012), Gautier and Tsybakov (2011)
allow to handle non-identically distributed errors $u_i$ and are {\it pivotal}, i.e., rely on very weak distributional assumptions.
In Gautier and Tsybakov (2011), the regressors $x_i$ can be correlated with
the errors $u_i$,  
and an estimator is suggested that
makes use of instrumental variables, called the {\em STIV} (Self-Tuned Instrumental Variables) estimator.
In a particular instance, the {\em STIV} estimator
can be applied in classical linear regression model where
all regressors are uncorrelated with the errors. This yields a pivotal extension
of the Dantzig selector based on conic programming.
Gautier and Tsybakov (2011) also present a method to obtain finite sample confidence sets
that are robust
to non-Gaussian and heteroscedastic errors. 

Another important issue is to relax the assumptions on the model under which the validity of the Lasso type methods is proved, such as the restricted eigenvalue condition of Bickel,
Ritov and Tsybakov (2009) and its various analogs. Belloni, Chernozhukov and Wang (2011b) obtain
fast rates for prediction for the Square-root Lasso under a relaxed version of the restricted eigenvalue
condition.  In the context of known noise variance, Ye and Zhang (2011) introduce cone invertibility factors instead of restricted eigenvalues. For pivotal estimation, an approach based on the sensitivities and sparsity certificates is introduced in Gautier and Tsybakov (2011), see more details below. Finally, note that aggregation by exponential weighting (Dalalyan and Tsybakov~(2008), Rigollet and
Tsybakov~(2011, 2012)) does not require any condition on the model but its numerical realization is based on MCMC algorithms in high dimension whose convergence rate is hard to assess theoretically.

In this paper, we introduce a new pivotal estimator, called the Self-tuned Dantzig estimator. It is defined as a linear program, so from the numerical point of view it is simpler than the previously known pivotal estimators based on conic programming. We obtain upper bounds on its estimation and prediction errors under weak assumptions on the model and on the distribution of the errors showing that it achieves the same rate as in the more restrictive situation of fixed design and i.i.d.
Gaussian errors with known variance. The model assumptions are based on the sensitivity analysis from Gautier and Tsybakov (2011). Distributional assumptions allow for dependence between $x_{i}$ and $u_i$. When $x_{i}$'s  are independent from $u_i$'s, it is enough to assume, for example, that the  errors $u_i$ are symmetric and have a finite second moment.

\section{Notation}
\label{s2}

We set $\bold{Y}=(y_1,\dots,y_n)^T$, $\bold{U}= (u_1,\dots,u_n)^T$,
and we denote by $\bold{X}$ the matrix of dimension
$n\times p$ with rows $x_i^T$, $i=1,\hdots,n$. We denote by ${\bf D}$ the  $p\times p$ diagonal normalizing matrix
with diagonal entries $d_{kk}>0$, $k=1,\dots,p$. Typical examples are: $d_{kk}\equiv 1$ or
$$
d_{kk} = 
\left(\frac1{n}\sum_{i=1}^n x_{ki}^2\right)^{-1/2}, \  \  \text{and} \  \ d_{kk} =
\Big(\max_{i=1,\dots,n} | x_{ki}|\Big)^{-1}
$$
where $x_{ki}$ is the $k$th component of $x_i$.
For a vector $\beta\in\R^p$, let $J(\beta)=\{k\in\{1,\hdots,p\}:\
\beta_k\ne0\}$ be its support, {i.e.}, the set of indices
corresponding to its non-zero components $\beta_k$. We denote by
$|J|$ the cardinality of a set $J\subseteq \{1,\hdots,p\}$ and by
$J^c$ its complement: $J^c=\{1,\hdots,p\}\setminus J$.
The $\ell_p$ norm of a vector $\Delta$ is denoted by
$|\Delta|_p$, $1\le p\le\infty$. For
$\Delta=(\Delta_1,\dots\Delta_p)^T\in\R^p$ and a set of indices
$J\subseteq \{1,\ldots,p\}$, we consider $\Delta_J\triangleq
(\Delta_1\indic_{\{1\in J\}}, \ldots,\Delta_p\indic_{\{p\in
J\}})^T$, where $\indic_{\{\cdot\}}$ is the indicator function.
For $a\in \R$, we set $a_+\triangleq\max(0,a)$, $a_+^{-1}\triangleq
(a_+)^{-1}$.

\section{The Estimator}\label{s3}
We say that a pair $(\beta,\sigma)\in\R^p\times\R^+$
satisfies the {\em Self-tuned Dantzig}-{\em constraint} if it belongs to the set
\begin{equation}\label{STV}
\widehat{\mathcal{D}}\triangleq\left\{(\beta,\sigma)\
\beta\in\R^p,\ \sigma>0,\ \left|\frac1n {\bf D}\bold{X}^T(\bold{Y}
-\bold{X}\beta)\right|_{\infty}\le \sigma
r\right\}
\end{equation}
for some $r>0$ (specified below).

\begin{dfntn}
We call the {\em Self-Tuned Dantzig} estimator any solution
$(\widehat{\beta},\widehat{\sigma})$ of the following minimization
problem
\begin{equation}\label{min1}
\min_{(\beta,\sigma)\in \widehat{\mathcal{D}}} \left(\,\left|{\bf
D}^{-1}\beta\right|_1+c\sigma \right),
\end{equation}
for some positive constant $c$.
\end{dfntn}

Finding the {Self-Tuned Dantzig} estimator is a linear
program. The term $c\sigma$ is included in the criterion
to prevent from choosing $\sigma$ arbitrarily large. The choice of the constant
$c$ will be discussed later.

\section{Sensitivity Characteristics}\label{s42}
The sensitivity characteristics are defined by the action
of the matrix
$$\Psi_n\triangleq\frac1n {\bf D}\bold{X}^T\bold{X}{\bf D}$$
on the so-called {\em cone of dominant coordinates}
$$
C_{J}^{(\gamma)}\triangleq\left\{\Delta\in\R^p:\
|\Delta_{J^c}|_1\le(1+\gamma)|\Delta_{J}|_1 \right\},
$$
for some $\gamma>0$.   It is straightforward that
for $\delta\in C_{J}^{(\gamma)}$,
\begin{equation}\label{eq:technical}
|\Delta|_1\le(2+\gamma)|\Delta_{J}|_1\le (2+\gamma)|J|^{1-1/q}
|\Delta_{J}|_q,\quad \forall 1\le q\le\infty.
\end{equation}
We now recall some definitions from
Gautier and Tsybakov (2011).
For $q\in[1,\infty]$, we define the {\em $\ell_q$ sensitivity} as
the following random variable
$$\kappa_{q,J}^{(\gamma)}\triangleq\inf_{\Delta\in C_{J}^{(\gamma)}:\
|\Delta|_q=1}\left|\Psi_n\Delta \right|_{\infty}.$$
Given a subset $J_0\subset\{1,\dots,p\}$ and $q\in[1,\infty]$,
we define the $\ell_q$-$J_0$-{\em block sensitivity} as
\begin{equation}\label{eq:technical2}
\kappa_{q,J_0,J}^{(\gamma)}\triangleq\inf_{\Delta\in C_J^{(\gamma)}:\
|\Delta_{J_0}|_q=1}\left|\Psi_n\Delta \right|_{\infty}.
\end{equation}
By convention, we set $\kappa_{q,{\varnothing} ,J}^{(\gamma)}=\infty$.
%
Also, recall that the restricted eigenvalue
of Bickel, Ritov and Tsybakov (2009) is defined by
$$
\kappa_{{\rm RE},J}^{(\gamma)} \triangleq \inf_{\Delta\in\R^p\setminus\{0\}:\
\Delta\in C_{J}^{(\gamma)}} \frac{|\Delta^T\Psi_n\Delta|}{|\Delta_{J}|_2^2}$$
and a closely related quantity is
$$\kappa^{'(\gamma)}_{{\rm RE},J} \triangleq \inf_{\Delta\in\R^p\setminus\{0\}:\
\Delta\in C_{J}^{(\gamma)}}
\frac{|J|\,|\Delta^T\Psi_n\Delta|}{|\Delta_{J}|_1^2}.
$$
The next result establishes a relation between restricted eigenvalues and sensitivities.  It follows directly from the Cauchy-Schwarz inequality and \eqref{eq:technical}.
\begin{lmm}\label{lemma1}
\begin{equation}\label{eq:compRE}
\kappa_{{\rm RE},J}^{(\gamma)}\le \kappa^{'(\gamma)}_{{\rm RE},J} \le (2+\gamma)|J|
\kappa_{1,J,J}^{(\gamma)}\le (2+\gamma)^2|J|\kappa_{1,J}^{(\gamma)}.
\end{equation}
\end{lmm}
The following proposition gives a useful lower bound
on the sensitivity. 
\begin{prpstn}\label{p1} If $|J|\le s$,
\begin{align}
\kappa_{1,J,J}^{(\gamma)}&\ge \frac{1}{s}\min_{k=1,\dots,p}\left\{
\min_{\ \ \Delta_k= 1,\
|\Delta|_1\le (2+\gamma)s}\left|\Psi_n\Delta\right|_{\infty}\right\} \triangleq
\kappa_{1,0}^{(\gamma)}(s)\label{eq:lb3}.
\end{align}
\end{prpstn}
{\bf Proof.} We have
\begin{align*}
\kappa_{1,J,J}^{(\gamma)}&=\inf_{\Delta:\ |\Delta_J|_1=1,\ |\Delta_{J^c}|_1
\le 1+\gamma}\left|\Psi_n\Delta \right|_{\infty}\\
&\ge \inf_{\Delta:\ |\Delta|_{\infty}\ge\frac{1}{s},\ |\Delta|_1\le 2+\gamma}
\left|\Psi_n\Delta \right|_{\infty}\\
&= \frac{1}{s}\inf_{\Delta:\ |\Delta|_{\infty}\ge1,\ |\Delta|_1\le (2+\gamma)s}
\left|\Psi_n\Delta \right|_{\infty}\quad ({\rm by\ homogeneity})\\
&= \frac{1}{s}\inf_{\Delta:\ |\Delta|_{\infty}\ge1,\ |\Delta|_1\le (2+\gamma)s}
|\Delta|_{\infty}\frac{\left|\Psi_n\Delta \right|_{\infty}}{|\Delta|_{\infty}}\\
&\ge \frac{1}{s}\inf_{\Delta:\ |\Delta|_{\infty}=1,\ |\Delta|_1\le (2+\gamma)s
|\Delta|_{\infty}}\left|\Psi_n\Delta \right|_{\infty}\quad ({\rm by\ homogeneity})\\
&= \frac{1}{s}\inf_{\Delta:\ |\Delta|_{\infty}=1,\ |\Delta|_1\le (2+\gamma)s}
\left|\Psi_n\Delta \right|_{\infty}\\
&= \frac{1}{s}\min_{k=1,\hdots,p}\left\{\inf_{\Delta:\ \Delta_k=1,\ |\Delta|_1
\le (2+\gamma)s}\left|\Psi_n\Delta \right|_{\infty}\right\}\,.
\quad\quad\quad\quad\quad\quad\quad\quad\square
\end{align*}
Note that the random variable $\kappa_{1,0}^{(\gamma)}(s)$ 
depends only on the observed data. It is not difficult to see that it can be obtained by solving $p$ linear programs.   For
more details and further results on the sensitivity characteristics, see Gautier and Tsybakov (2011).

\section{Bounds on the estimation and prediction errors}\label{s5}
In this section, we use the notation $\Delta\triangleq {\bf D}^{-1}(\widehat{\beta} -\beta)$.
Let $0<\alpha<1$ be a given constant. 
We choose
the tuning parameter $r$ in the definition of $\widehat{\mathcal D}$ as follows:
\begin{equation}\label{er1}
r=\sqrt{\frac{2\log(4p/\alpha)}{n}}.
\end{equation}

\begin{thrm}\label{l11}Let for all $i=1,\hdots,n$, and $k=1,\hdots,p$,
the random variables $x_{ki}u_i$ be
symmetric. 
Let $Q^*>0$ be a constant such that
\begin{equation}\label{eq:Qcond}
\mathbb{P}\left(\max_{k=1,\hdots,p}  \frac{d_{kk}^2}{n}\sum_{i=1}^n{x_{ki}^2}
u_i^2>Q^*\right)\le \alpha/2.
\end{equation}
 Assume that $|J(\beta^*)|\le s$, and set in (\ref{min1})
\begin{equation}\label{eq:cvalue0}
c=\frac{(2\gamma+1)r}{\kappa_{1,0}^{(\gamma)}(s)},
\end{equation}
where $\gamma$ is a positive number. Then, with probability at least   $1-\alpha$, for any $\gamma>0$ and
any $\widehat\beta$ such that $(\widehat\beta,\widehat\sigma)$
is a solution of the minimization problem \eqref{min1} with $c$ defined in (\ref{eq:cvalue0}) we have the following bounds on the $\ell_1$ estimation error and on the prediction error:
\begin{align}
\left|\Delta\right|_1&\le\left(\frac{(\gamma+2)(2\gamma+1)\sqrt{Q^*}}
{\gamma\kappa_{1,0}^{(\gamma)}(s)}\right) r\,, \label{eq:ub:estim}\\
\Delta^T\Psi_n\Delta &\le
\left(\frac{(\gamma+2)(2\gamma+1)^2Q^* }
{\gamma^2\kappa_{1,0}^{(\gamma)}(s)}\right)r^2\label{eq:ub:predict}.
\end{align}
\end{thrm}

\noindent{\bf Proof.} Set
\begin{align*}
\widehat{Q}(\beta)&\triangleq\max_{k=1,\hdots,p} \frac{d_{kk}^2}{n}\sum_{i=1}^n  x_{ki}^2
(y_i-x_i^T\beta)^2,
\end{align*}
and define the event
$$\mathcal{G}=\left\{\left|\frac1n{{\bf D}}{\bf X}^T
{\bf U}\right|_{\infty}\le
r\sqrt{\widehat{Q}(\beta^*)}\right\}=\left\{\left|\frac{d_{kk}}{n}\sum_{i=1}^nx_{ki}
u_i\right|\le
r\sqrt{\widehat{Q}(\beta^*)}, \ \ k=1,\hdots,p
\right\}\,.$$
Then
\begin{equation*}
\mathcal{G}^c\subset\bigcup_{k=1,\hdots,p}\left\{
\left|\frac{\sum_{i=1}^nx_{ki}u_i}{\sqrt{\sum_{i=1}^n(x_{ki}u_i)^2}}\right|\ge
\sqrt{n}r\right\}\,
\end{equation*}
and the union bound yields
\begin{equation}\label{eq:union}
\xP(\mathcal{G}^c)\le\sum_{k=1}^p\xP\left(
\left|\frac{\sum_{i=1}^nx_{ki}u_i}{\sqrt{\sum_{i=1}^n(x_{ki}u_i)^2}}\right|\ge
\sqrt{n}r\right).
\end{equation}
We now use the following result on deviations of self-normalized sums
due to Efron~(1969).
\begin{lmm}\label{th:efron} If $\eta_1,\dots, \eta_n$ are independent symmetric random
variables, then  
$$\mathbb{P}\left(\frac{\left|\frac{1}{n}\sum_{i=1}^n\eta_i\right|}
{\sqrt{\frac{1}{n}\sum_{i=1}^n\eta_i^2}}\ge t\right)\le
2\exp\left(-\frac{nt^2}{2}\right), \ \forall \ t>0.$$
\end{lmm}
For each of the probabilities on the right-hand side of (\ref{eq:union}), we apply Lemma~\ref{th:efron} with $\eta_i=x_{ki}u_i$. This and the definition of $r$ yield
$\mathbb{P}(\mathcal{G}^c)\le \alpha/2$. Thus, the event $\mathcal{G}$
holds with probability at least $1-\alpha/2$.
On the event $\mathcal{G}$ we have
\begin{align}
\left|\Psi_n\Delta \right|_{\infty} &\le \left|\frac{1}{n}{{\bf D}}
\bold{X}^T(\bold{Y}-\bold{X}\widehat{\beta})\right|_{\infty}
+\left|\frac{1}{n}{{\bf D}}\bold{X}^T(\bold{Y}-\bold{X}\beta^*)
\right|_{\infty}\label{stop1}\\
&\le r\widehat{\sigma}+\left|\frac1n{{\bf D}}{\bf X}^T
{\bf U}\right|_{\infty}\label{eq:Notobvious2}\\
&\le
r\left(\widehat{\sigma}+\sqrt{
\widehat{Q}(\beta^*)}\right)\nonumber\\
&\le
r\left[2\sqrt{
\widehat{Q}(\beta^*)}+\left(\widehat{\sigma}-\sqrt{
\widehat{Q}(\beta^*)}\right)\right]\label{stop2}
\end{align}
Inequality \eqref{eq:Notobvious2} holds because
$(\widehat{\beta},\widehat{\sigma})$ belongs to the
set $\widehat{\mathcal{D}}$ by definition.  Notice that, on
the event $\mathcal{G}$,
$\left(\beta^*,\sqrt{\widehat{Q}(\beta^*)}\right)$ belongs to the
set $\widehat{\mathcal{D}}$.  On the other hand,
$(\widehat{\beta},\widehat{\sigma})$ minimizes the
criterion $\left|{{\bf D}}^{-1}\beta\right|_1 +c\sigma$ on the same
set $\widehat{\mathcal{D}}$. Thus, on the event $\mathcal{G}$,
\begin{equation}\label{eq:main}
\left|{{\bf D}}^{-1}\widehat{\beta}\right|_1 +c\widehat{\sigma}\le
|{{\bf D}}^{-1}\beta^*|_1+c\sqrt{\widehat{Q}(\beta^*)}.
\end{equation}
This implies, again on the event $\mathcal{G}$,
\begin{align}
\left|\Psi_n\Delta \right|_{\infty}&\le
r\left[2\sqrt{
\widehat{Q}(\beta^*)}+\frac{1}{c}
\sum_{k\in J(\beta^*)}\left(
\left|d_{kk}^{-1}\beta^*_k\right|
-\left|d_{kk}^{-1}\widehat{\beta}_k\right|
\right)-\frac{1}{c}\sum_{k\in J(\beta^*)^c}\left|d_{kk}^{-1}
\widehat{\beta}_k\right|\right]\nonumber\\
&\le
r\left(2\sqrt{
\widehat{Q}(\beta^*)}+\frac{1}{c}\left|\Delta_{J(\beta^*)}\right|_1\right) \label{stop2b}
\end{align}
where $\beta^*_k, \widehat{\beta}_k$ are the $k$th components of $\beta^*, \widehat{\beta}$. Similarly, (\ref{eq:main}) implies that, on the event $\mathcal{G}$,
\begin{align}
\left|\Delta_{J(\beta^*)^c}\right|_1&=\sum_{k\in J(\beta^*)^c}
\left|d_{kk}^{-1}\widehat{\beta}_k\right|\nonumber\\
&\le\sum_{k\in J(\beta^*)}\left(
\left|d_{kk}^{-1}\beta^*_k\right|
-\left|d_{kk}^{-1}\widehat{\beta}_k\right|\right)
+c\left(\sqrt{\widehat{Q}(\beta^*)}-
\widehat{\sigma}\right)\nonumber\\
&\le \left|\Delta_{J(\beta^*)}\right|_1+c\sqrt{\widehat{Q}(\beta^*)}\label{in0}.
\end{align}
We now distinguish between the following two cases.\\
Case 1: $c\sqrt{\widehat{Q}(\beta^*)}\le\gamma\left|\Delta_{J(\beta^*)}\right|_1$. In this case
\eqref{in0} implies
\begin{equation}\label{eq:conec}
\left|\Delta_{J(\beta^*)^c}\right|_1\le
(1+\gamma)\left|\Delta_{J(\beta^*)}\right|_1.
\end{equation}
Thus, $\Delta\in C_{J(\beta^*)}^{(\gamma)}$ on the event $\mathcal{G}$.
By definition of $\kappa_{1,J(\beta^*)
,J(\beta^*)}^{(\gamma)}$ and (\ref{eq:lb3}),
$$\left|\Delta_{J(\beta^*)}\right|_1\le\frac{\left|\Psi_n\Delta \right|_{\infty}}{\kappa_{1,J(\beta^*)
,J(\beta^*)}^{(\gamma)}}\le \frac{\left|\Psi_n\Delta \right|_{\infty}}{
\kappa_{1,0}^{(\gamma)}(s)}\,.$$
This and \eqref{stop2b} yield
$$\left|\Delta_{J(\beta^*)}\right|_1\le\frac{2r\sqrt{\widehat{Q}(\beta^*)}}
{\kappa_{1,0}^{(\gamma)}(s)}
\left(1-\frac{r}{c\kappa_{1,0}^{(\gamma)}(s)}\right)_+^{-1}\,.$$
Case 2: $c\sqrt{\widehat{Q}(\beta^*)}>\gamma\left|\Delta_{J(\beta^*)}\right|_1$.
Then, obviously,
$\left|\Delta_{J(\beta^*)}\right|_1<\frac{c}{\gamma}\sqrt{\widehat{Q}(\beta^*)}.$

Combining the two cases we obtain, on the event $\mathcal{G}$,
\begin{equation}\label{eq:ubDeltaJ}
\left|\Delta_{J(\beta^*)}\right|_1\le\sqrt{\widehat{Q}(\beta^*)}\max\left\{\frac{2r}
{\kappa_{1,0}^{(\gamma)}(s)}
\left(1-\frac{r}{c\kappa_{1,0}^{(\gamma)}(s)}\right)_+^{-1}, \,
\frac{c}{\gamma}\right\}\,.
\end{equation}
In this argument, $c>0$ and $\gamma>0$ were arbitrary. The value of $c$ given in (\ref{eq:cvalue0}) is the minimizer of the
right-hand side of (\ref{eq:ubDeltaJ}).  Plugging it in (\ref{eq:ubDeltaJ}) we find that, with probability at least $1-\alpha/2$
\begin{align*}
\left|\Delta\right|_1&\le\frac{(\gamma+2)(2\gamma+1)r}
{\gamma\kappa_{1,0}^{(\gamma)}(s)}\sqrt{\widehat{Q}(\beta^*)}
\end{align*}
where we have used \eqref{in0}.  Now, by (\ref{eq:Qcond}),  $\widehat{Q}(\beta^*)\le Q^*$ with probability at least $1-\alpha/2$. Thus, we get that (\ref{eq:ub:estim}) holds with probability at least $1-\alpha$. Next, using \eqref{stop2b} we obtain that, on the same event of probability at least $1-\alpha$,
\begin{align*}
\left|\Psi_n\Delta \right|_{\infty}&\le
\frac{(2\gamma+1)r}{\gamma}\sqrt{Q^*}.
\end{align*}
Combining this inequality with (\ref{eq:ub:estim}) yields (\ref{eq:ub:predict}).
 \hfill$\square$\vspace{0.3cm}

{\it Discussion of Theorem~\ref{l11}.}
\begin{enumerate}
\item In view of Lemma~\ref{lemma1}, $\kappa_{1,J(\beta^*)
,J(\beta^*)}^{(\gamma)}\ge (2+\gamma)^{-2} \kappa_{{\rm RE},J(\beta^*)}^{(\gamma)}/s$. Also, it is easy to see from Proposition~\ref{p1} that $\kappa_{1,0}^{(\gamma)}(s)$ is of the order $1/s$ when $\Psi_n$ is the identity matrix and $p\gg s$ (this is preserved for $\Psi_n$ that are small perturbations of the identity). Thus, 
the bounds (\ref{eq:ub:estim}) and (\ref{eq:ub:predict})  take the form
$$
\left|\Delta\right|_1 \le C\left(s\sqrt{\frac{\log p}{n}}\right), \qquad
\Delta^T\Psi_n\Delta \le C\left({\frac{s\log p}{n}}\right)\,,
$$
for some constant $C$, and we recover the usual rates for the $\ell_1$ estimation and for the prediction error respectively, cf. Bickel, Ritov and Tsybakov (2009).
\item Theorem~\ref{l11} does not assume that $x_{ki}$'s are independent from $u_i$'s. The only assumption is the symmetry of $x_{ki}u_i$. However, if $x_{ki}$ is independent from $u_i$, then by conditioning on $x_{ki}$ in the bound for ${\mathbb P}({\mathcal G})$, it is enough to assume the symmetry of $u_i$'s. Furthermore, while we have chosen the symmetry since it makes the conditions of Theorem~\ref{l11} simple and transparent, it is not essential for our argument to be applied. The only point in the proof where we use the symmetry is the bound for the probability of deviations of self-normaized sums  ${\mathbb P}({\mathcal G})$. This probability can be bounded in many other ways without the symmetry assumption, cf., e.g., Gautier and Tsybakov~(2011). It is enough to have $\E[x_{ki}u_i]=0$ and a uniform over $k$ control of the ratio
$$\frac{(\sum_{i=1}^n
\E[x_{ki}^2u_i^2])^{1/2}}
{\left(\sum_{i=1}^n\E[|x_{ki}u_{i}|^{2+\delta}]\right)^{1/(2+\delta)}}
$$
for some $\delta>0$, cf. \cite{JSW} or \cite{BGHK2}.
\item The quantity $Q^*$ is not present in the definition of the estimator and is needed only to assess the rate of convergence. It is not hard to find $Q^*$ in various situations.  The simplest case is when $d_{kk}\equiv 1$ and the random variables $x_{ki}$ and $ u_i$ are bounded uniformly in $k,i$ by a constant $L$. Then we can take $Q^*=L^4$. If only $x_{ki}$ are bounded uniformly in $k$ by $L$, condition (\ref{eq:Qcond}) holds when
$\mathbb{P}\left(\frac{1}{n}\sum_{i=1}^n
u_i^2>Q^*/L^2\right)\le \alpha/2,$ and then for $Q^*$ to be bounded it is enough to assume that $u_i$'s have a finite second moment. The same remark applies when $d_{kk} = \left(\max_{i=1,\dots,n} | x_{ki}|\right)^{-1}$, with an advantage that in this case we guarantee that $Q^*$ is bounded under no assumption on $ x_{ki}$.
\item The bounds in Theorem~\ref{l11} depend on $\gamma>0$ that can be optimized. Indeed, the functions of $\gamma$ on the right-hand sides of (\ref{eq:ub:estim}) and (\ref{eq:ub:predict}) are data-driven and can be minimized on a grid of values of $\gamma$. Thus, we obtain an optimal (random) value $\gamma=\hat\gamma$, for which (\ref{eq:ub:estim}) and (\ref{eq:ub:predict}) remain valid, since these results hold for any $\gamma>0$.
\end{enumerate}





\begin{thebibliography}{99.}%
\bibitem{Ant}
{Antoniadis, A.} (2010)
Comments on: $l^1$-penalization for Mixture Regression Models.
(with discussion). {\em Test,} {\bf 19}, 257--258.
\bibitem{BCCH}
{Belloni, A., Chen, D., Chernozhukov, V. and Hansen C.} (2012)
 Sparse Models and Methods for Optimal Instruments
with an Application to Eminent Domain. {\em Econometrica,}
{\bf 80}, 2369--2430.
\bibitem{BC3}
{Belloni, A. and Chernozhukov, V.} (2010) Least
Squares After Model Selection in High-dimensional Sparse
Models. Forthcoming in {\em Bernoulli}.
\bibitem{BCWa}
{Belloni, A., Chernozhukov, V. and Wang, L.} (2011a)
Square-Root Lasso: Pivotal Recovery of Sparse
Signals via Conic Programming. {\em Biometrika,}
{\bf 98}, 791--806.
\bibitem{BCWb}
{Belloni, A., Chernozhukov, V. and Wang, L.}  (2011b)
 Pivotal Estimation of Nonparametric Functions via
Square-root Lasso.
Preprint \url{http://arxiv.org/pdf/1105.1475.pdf} 
\bibitem{BGHK2} Bertail, P. , Gauth\'erat, E. and
Harari-Kermadec, H. (2009) Exponential Inequalities
for Self Normalized Sums. \emph{Electronic
Communications in Probability,} {\bf 13}, 628--640.
\bibitem{BRT}
{Bickel, P.,  Ritov, J. Y. and Tsybakov, A. B. } (2009)
 Simultaneous Analysis of Lasso and Dantzig
Selector. {\em The Annals of Statistics,}  {\bf 37},
1705--1732.
\bibitem{BvdG}
{B\"uhlmann, P. and van de Geer, S.A.} (2011) {\em Statistics for
High-Dimensional Data}. Springer, New-York.
\bibitem{CT}
{Cand\`es, E., and Tao, T.} (2007)  The Dantzig
Selector: Statistical Estimation when $p$ is Much Larger than
$n$. {\em The Annals of Statistics,} {\bf 35},
2313--2351.
\bibitem{D}
{Dalalyan, A.} (2012)
SOCP Based Variance Free Dantzig Selector with Application to Robust Estimation.
{\em C. R. Math. Acad. Sci. Paris,} {\bf 350}, 785--788
\bibitem{DT}
{Dalalyan, A., and Tsybakov, A.B.} (2008)
Aggregation by Exponential Weighting, Sharp PAC-Bayesian Bounds and
Sparsity. {\em Journal of Machine Learning
Research,} {\bf 72}, 39--61.
\bibitem{Efron}
{Efron, B.} (1969) Student's t-test Under
Symmetry Conditions. {\em Journal of American
Statistical Association,} {\bf  64}, 1278--1302.
\bibitem{GT}
{Gautier, E. and Tsybakov, A.B.} (2011)
 High-dimensional instrumental
variables regression and confidence sets.
Preprint \url{http://arxiv.org/pdf/1105.2454v1.pdf}  
\bibitem{JSW}
{Jing, B.-Y.,  Shao, Q. M. and Wang, Q. } (2003)
Self-Normalized Cram\'er-Type Large Deviations for Independent
Random Variables. {\em The Annals of Probability,}
{\bf  31}, 2167--2215.
\bibitem{Kol1}
{Koltchinskii, V.} (2011)   {\em Oracle Inequalities for Empirical Risk Minimization and Sparse Recovery Problems}. Lecture Notes in Mathematics, vol. 2033. Springer, New-York.
\bibitem{RiT}
{Rigollet, P. and Tsybakov, A.B.} (2011)
Exponential Screening and Optimal Rates of Sparse
Estimation. \emph{The Annals of Statistics,} {\bf 39},
731--771.
\bibitem{RiT2}
{Rigollet, P. and Tsybakov, A.B.} (2012)
Sparse Estimation by Exponential Weighting.
{\it Statistical Science}, {\bf 27}, 558--575.
\bibitem{SBvdG}
{St\"adler, N., B\"uhlmann, P. and van de Geer, S.A.} (2010)
 $l^1$-penalization for Mixture Regression
Models. {\em Test,} {\bf 19}, 209--256.
\bibitem{SunZhang}
{Sun, T. and Zhang, C.-H.} (2011)
Scaled Sparse Linear Regression.
Preprint \url{http://arxiv.org/abs/1104.4595}
\bibitem{YZ} {Ye, F.
and Zhang, C.-H. } (2010) Rate Minimaxity of the
Lasso and Dantzig Selector for the $l_q$ Loss in $l_r$
Balls. {\em Journal of Machine Learning Research,}
{\bf 11}, 3519--3540.


\end{thebibliography}
\end{document}